\documentclass[11pt]{amsart}
\usepackage{mathrsfs}
\usepackage{amssymb,latexsym}
\usepackage{pstcol,pstricks,color}

 \numberwithin{equation}{section}

\newtheorem{theorem}{Theorem}[section]

\newtheorem{corollary}[theorem]{Corollary}

\newtheorem{remark}[theorem]{Remark}

\newtheorem{lemma}[theorem]{Lemma}

\def\qed{\hfill $\Box$}
\def\pf{\noindent {\it Proof.} }

\title{ Minors of a Class of Riordan Arrays Related to Weighted Partial Motzkin Paths }

\begin{document}
\maketitle
\begin{center}
Yidong Sun$^\dag$ and Luping Ma$^\ddag$


Department of Mathematics, Dalian Maritime University, 116026 Dalian, P.R. China\\[5pt]

{\it  Emails: $^\dag$sydmath@yahoo.com.cn }

\end{center}\vskip0.2cm

\subsection*{Abstract} A partial Motzkin path is a path from $(0, 0)$ to $(n, k)$ in the
$XOY$-plane that does not go below the $X$-axis and consists of up
steps $U=(1, 1)$, down steps $D=(1, -1)$ and horizontal steps $H=(1,
0)$. A weighted partial Motzkin path is a partial Motzkin path with
the weight assignment that all up steps and down steps are weighted
by $1$, the horizontal steps are endowed with a weight $x$ if they
are lying on $X$-axis, and endowed with a weight $y$ if they are not
lying on $X$-axis. Denote by $M_{n,k}(x, y)$ to be the weight
function of all weighted partial Motzkin paths from $(0, 0)$ to $(n,
k)$, and $\mathcal{M}=(M_{n,k}(x,y))_{n\geq k\geq 0}$ to be the
infinite lower triangular matrices. In this paper, we consider the
sums of minors of second order of the matrix $\mathcal{M}$, and
obtain a lot of interesting determinant identities related to
$\mathcal{M}$, which are proved by bijections using weighted partial
Motzkin paths. When the weight parameters $(x, y)$ are specialized,
several new identities are obtained related to some classical
sequences involving Catalan numbers. Besides, in the alternating
cases we also give some new explicit formulas for Catalan numbers.

\medskip

{\bf Keywords}: Motzkin path; Catalan number; Motzkin number;
Riordan array.

\noindent {\sc 2000 Mathematics Subject Classification}: Primary
05A19; Secondary 05A15, 05A10.

{\bf \section{ Introduction } }

The starting point for this paper is the observation that the close
connections between the Catalan numbers
$C_n=\frac{1}{n+1}\binom{2n}{n}$ and the Pascal triangle
$\mathscr{P}=(\binom{n}{k})_{n\geq k\geq 0}$, that is,
\begin{eqnarray*}
C_{n+1}=\sum_{k=0}^{n}N_{n+1,k+1} \hskip-.22cm &=&\hskip-.22cm
\sum_{k=0}^{n}\det\left(\begin{array}{cc}
\binom{n}{k} & \binom{n}{k+1} \\[5pt]
\binom{n+1}{k} & \binom{n+1}{k+1}
\end{array}\right), \\[5pt]
C_{n+1}=\sum_{k=0}^{n}N_{n+1,k+1} \hskip-.22cm &=&\hskip-.22cm
\sum_{k=0}^{n}\det\left(\begin{array}{cc}
\binom{n}{k} & \binom{n+1}{k+1} \\[5pt]
\binom{n+1}{k} & \binom{n+2}{k+1}
\end{array}\right),
\end{eqnarray*}
where $N_{n,k}=\frac{1}{n}\binom{n}{k}\binom{n}{k-1}$ for $n\geq
k\geq 1$ are the Narayana numbers.

This makes us to do some numerical verification for other combinatorial triangles. For example, consider
Shapiro's Catalan triangle \cite{ShapA}, defined by $\mathcal{B}=(B_{n,k})_{n\geq k\geq 0}$ such that
$B_{n,k}=\frac{k+1}{n+1}\binom{2n+2}{n-k}$. Table 1 illustrates this triangle for small $n$ and $k$ up to $5$.
\begin{center}
\begin{eqnarray*}
\begin{array}{c|cccccccc}\hline
n/k & 0   & 1    & 2    & 3    & 4    & 5       \\\hline
  0 & 1   &      &      &      &      &        \\
  1 & 2   & 1    &      &      &      &         \\
  2 & 5   & 4    & 1    &      &      &         \\
  3 & 14  & 14   & 6    &  1   &      &         \\
  4 & 42  & 48   & 27   &  8   &  1   &        \\
  5 & 132 & 165  & 110  &  44  &  10  &  1     \\\hline
\end{array}
\end{eqnarray*}
Table 1. The values of $B_{n,k}$ for $n$ and $k$ up to $5$.
\end{center}

Let $\mathcal{X}=(X_{n,k})_{n\geq k\geq 0}$ be the infinite lower triangles defined on the triangle $\mathcal{B}$ by
\begin{eqnarray*}
X_{n,k} \hskip-.22cm &=&\hskip-.22cm  \rm{det}\left(\begin{array}{cc}
B_{n, k}   & B_{n,k+1}  \\[5pt]
B_{n+1,k}  & B_{n+1,k+1}
\end{array}\right).
\end{eqnarray*}
Table 1.2 illustrates the triangle $\mathcal{X}$ for small $n$ and $k$ up to $4$, together
with the row sums. It indicates that the row sums have close relation with the first column of the triangle $\mathcal{B}$.
\begin{center}
\begin{eqnarray*}
\begin{array}{c|ccccc|c}\hline
n/k & 0   & 1   & 2    & 3    & 4     & row\ sums           \\\hline
  0 & 1   &     &      &      &       &  1=1^2              \\
  1 & 3   & 1   &      &      &       &  4=2^2             \\
  2 & 14  & 10  & 1    &      &       &  25=5^2             \\
  3 & 84  & 90  & 21   & 1    &       &  196=14^2            \\
  4 & 594 & 825 & 308  & 36   &  1    &  1764=42^2           \\\hline
\end{array}
\end{eqnarray*}
Table 1.2. The values of $X_{n,k}$ for $n$ and $k$ up to $4$, together with the row sums.
\end{center}

The above fact motivates us to consider the following problem.

{\bf Question: } Let $\mathcal{A}=(A_{n,k})_{n\geq k\geq 0}$ be an
infinite lower triangular matrix with nonzero entries on the main
diagonal. Given integers $m, r, \ell, p$ with $m, \ell, p\geq 0$, define a transformation on $\mathcal{A}$ by
$\mathcal{A}_{p}=\big(A_{n,k}^{(p)}(m, r, \ell)\big)_{n\geq k\geq 0}$, where
\begin{eqnarray*}
A_{n,k}^{(p)}(m, r, \ell) \hskip-.22cm &=&\hskip-.22cm \det\left(A_{n+im+jr, k+j\ell}\right)_{0\leq i, j\leq p}.
\end{eqnarray*}
Then how to determine the explicit expression for the $n$-th row sum of $\mathcal{A}_{p}$,
$$S_{n,m,r,\ell}^{(p)}(\mathcal{A}_p)=\sum_{k=0}^{n}A_{n,k}^{(p)}(m, r, \ell)? $$

In general, it is not easy to give an exact answer for this
question. But, in the case $p=1$, for some special infinite lower
triangular matrices related to weighted partial Motzkin paths, it
can produce several surprising results.

The organization of this paper is as follows. The next section gives
a brief introduction to weighted partial Motzkin paths, which
generates a class of infinite lower triangular matrices. In Section
3, we state our main results and give bijective proofs. When the
weight parameters are specialized, several new identities are
obtained related to some classical sequences involving Catalan
numbers. In Section 4, we consider the alternating sums and give
some new explicit formulas for Catalan numbers. \vskip.2cm

\section{ Weighted partial Motzkin paths }

Recall that a {\em Motzkin path} is a lattice path from $(0, 0)$ to
$(n, 0)$ in the $XOY$-plane that does not go below the $X$-axis and
consists of up steps $U=(1, 1)$, down steps $D=(1, -1)$ and
horizontal steps $H=(1, 0)$. A {\em partial Motzkin path}, also
called a Motzkin path from $(0, 0)$ to $(n, k)$ in \cite{CamNk}, is
just a Motzkin path but without the requirement of ending on the
X-axis. A {\em weighted partial Motzkin path} \cite{Sloane} is a
partial Motzkin path with the weight assignment that the all up
steps and down steps are weighted by $1$, the horizontal steps are
endowed with a weight $x$ if they are lying on $X$-axis, and endowed
with a weight $y$ if they are not lying on $X$-axis. The {\em
weight} $w(P)$ of a path $P$ is the product of the weight of all its
steps. The {\em weight} of a set of paths is the sum of the total
weights of all the paths. If $P=L_1L_2\dots L_{n-1}L_{n}$ is
aweighted partial Motzkin path of length $n$, denoted by
$\overline{P}=\overline{L}_{n}\overline{L}_{n-1}\dots
\overline{L}_2\overline{L}_1$ the reverse of the path $P$, where
$\overline{L}_i=U$ if $L_i=D$, $\overline{L}_i=D$ if $L_i=U$ and
$\overline{L}_i=H$ if $L_i=H$. For any step, we say that it is at
level $i$ if the $Y$-coordinate of its end point is $i$. An up step
at level $i$ is $R$-$visible$ \cite{ChenLi} if it is the rightmost
up step at level $i$ and there are no other up steps at the same
level to its right. It is also worth mentioning that another type of
weighted partial Motzkin paths is used by Chen, Li, Shapiro and Yan
\cite{ChenLi} to derive many nice matrix identities related to a
class of Riordan arrays.

Let $\mathscr{M}_{n,k}(x,y)$ denote the set of weighted partial
Motzkin paths from $(0, 0)$ to $(n, k)$, and $M_{n,k}(x,y)$ be its
weight. Define $\mathscr{M}_{k}(x,y)=\bigcup_{n\geq
k}\mathscr{M}_{n,k}(x,y)$, it should be pointed out that any $P\in
\mathscr{M}_{k}(x,y)$ has exactly $k$ R-visible up steps. For any
$P\in \mathscr{M}_{n,k}(x,y)$, according to the last step $U, H$ or
$D$ of $P$, one can easily deduce the following recurrences for
$M_{n,k}(x,y)$,
\begin{eqnarray}
M_{n, 0}(x,y) \hskip-.25cm &=& \hskip-.25cm xM_{n-1,0}(x,y)+M_{n-1,1}(x,y), \ (n\geq 1),  \label{eqn 2.0}\\
M_{n, k}(x,y) \hskip-.25cm &=& \hskip-.25cm
M_{n-1,k-1}(x,y)+yM_{n-1,k}(x,y)+M_{n-1,k+1}(x,y),\ (n\geq k\geq 1),
\nonumber
\end{eqnarray}
with $M_{0,0}(x,y)=1$ and $M_{n, k}(x,y)=0$ if $n<k$ or $k<0$.

\begin{center}
\begin{eqnarray*}
\begin{array}{|c|l|l|l|l|l|}\hline
n/k & 0 & 1 & 2 & 3 & 4 \\\hline
0 & 1 & & & & \\
1 & x & 1 & & & \\
2 & x^2+1 & x+y & 1 & & \\
3 & x^3+2x+y & x^2+xy+y^2+2 & x+2y & 1 & \\
4 & x^4+3x^2+2xy+y^2+2 & x^3+x^2y+xy^2+3x+y^3+5y & x^2+2xy+3y^3+3 &
x+3y & 1 \\\hline
\end{array}
\end{eqnarray*}
Table 1. The values of $M_{n,k}(x,y)$ for $n$ and $k$ up to $4$.
\end{center}

Denote $(M_{n,k}(x,y))_{n\geq k\geq 0}$ by $\mathcal{M}$, then
$\mathcal{M}$ is an infinite lower triangular matrix with the main
diagonal entries $1$. Table 1 illustrates this matrix for small $n$
and $k$ up to $4$. In fact, the matrix $\mathcal{M}$ forms a Riordan
array. Recall that $\mathcal{R}=(R_{n,k})_{n\geq k\geq 0}$ is a
Riordan array \cite{ShapB, ShapGet, Sprug} if it is an infinite
lower triangular matrix with nonzero entries on the main diagonal,
such that $R_{n, k}=[t^n]g(t)(f(t))^{k}$ for $n\geq k$, namely,
$R_{n,k}$ equals the coefficient of $t^n$ in the expansion of the
series $g(t)(f(t))^k$, where $g(t)=1+g_1t +g_2t^2+\cdots $ and
$f(t)=f_1t + f_2t^2 + \cdots $ with $f_1\neq 0$ are two formal power
series. It is convenient to denote the Riordan array $\mathcal{R}$
by $(g(t), f(t))$. Let $M_k(x,y; t)=\sum_{n\geq k}M_{n, k}(x,y)t^n$
be the generating function of weighted partial Motzkin paths ending
at level $k$. For any $P\in \mathscr{M}_{k}(x,y)$, according to the
$k$ R-visible up steps for $k\geq 1$, $P$ can be uniquely
partitioned into $P=P_0UP_1U\dots UP_k$, where $P_0\in
\mathscr{M}_{0}(x,y)$ and $P_i\in \mathscr{M}_{0}(y,y)$ for $1\leq
i\leq k$. This decomposition produces a relation between $M_k(x,y;
t)$ and $M_0(x,y; t)$, namely,
\begin{eqnarray*}
M_{k}(x,y; t) \hskip-.25cm &=& \hskip-.25cm M_{0}(x,y;
t)(tM_{0}(y,y; t))^k, \ (k\geq 1),
\end{eqnarray*}
which indicates that $\mathcal{M}=(M_{0}(x,y; t), tM_{0}(y,y; t))$
is a Riordan array.

For any $P\in \mathscr{M}_{0}(x,y)$, $P$ has three cases to be
considered, that is (1) $P=\varepsilon$, an empty path; (2) starting
with a horizontal step, i.e., $P=HP_1$, where $P_1\in
\mathscr{M}_{0}(x,y)$; (3) $P=UP_2DP_1$, where $P_1\in
\mathscr{M}_{0}(x,y)$ and $P_2\in \mathscr{M}_{0}(y,y)$. Making use
of the so-called symbol method (for details see \cite{Sedg}), we
obtain
\begin{eqnarray*}
M_{0}(x,y; t) \hskip-.25cm &=& \hskip-.25cm 1+xtM_{0}(x,y;
t)+t^2M_{0}(y,y; t)M_{0}(x,y; t).
\end{eqnarray*}
Solving this equation, we have
\begin{eqnarray}\label{eqn 2.1}
M_{0}(y,y; t) \hskip-.25cm &=& \hskip-.25cm \frac{1-yt-\sqrt{(1-yt)^2-4t^2}}{2t^2}, \label{eqn 2.1} \\
M_{0}(x,y; t) \hskip-.25cm &=& \hskip-.25cm
\frac{1}{1-xt-t^2M_{0}(y,y; t)}
=\frac{1-2xt+yt-\sqrt{(1-yt)^2-4t^2}}{2(y-x)(1-xt)t+2t^2}.
\label{eqn 2.2}
\end{eqnarray}

When the parameters $x$ and $y$ are specialized, $M_{0}(x,y; t)$
produces generating functions for many classical combinatorial
sequences. We give a short list in Table 2, where
$C(t)=\frac{1-\sqrt{1-4t}}{2t}$ and
$M(t)=\frac{1-t-\sqrt{1-2t-3t^2}}{2t^2}$ are generating functions
respectively for Catalan numbers $C_n$ and Motzkin numbers $M_n$.
\begin{center}
\begin{eqnarray*}
\begin{array}{c|c|l|c}\hline
(x,y) & M_{0}(x,y; t)                        & Sequences                               & \mathcal{M} \\[2pt]\hline
(0,0) & C(t^2)                               & A126120=C_{\frac{n}{2}}                 & A053121 \  \mbox{\cite{Sloane} } \\[2pt]\hline
(0,1) & \frac{1+t-\sqrt{1-2t-3t^2}}{2t(1+t)} & A005043=Riordan\ numbers \ R_n          & A089942 \  \mbox{\cite{Sloane} } \\[2pt]\hline
(0,2) & \frac{1-\sqrt{1-4t}}{3-\sqrt{1-4t}}  & A000957=Fine\ numbers \ F_n             & A126093 \  \mbox{\cite{Sloane} } \\[2pt]\hline
(1,0) & \frac{1}{1-t-t^2C(t^2) }             & A001405=\binom{n}{\lfloor n/2\rfloor }  & A061554 \  \mbox{\cite{Sloane} } \\[2pt]\hline
(1,1) & M(t)                                 & A001006=Motzkin\ numbers\ M_n           & A064189 \  \mbox{\cite{Sloane} } \\[2pt]\hline
(1,2) & C(t)                                 & A000108=Catalan\ numbers\ C_n           & A039599 \  \mbox{\cite{Sloane} } \\[2pt]\hline
(2,2) & C^2(t)                               & A000108=Catalan\ numbers\ C_{n+1}       & A039598 \  \mbox{\cite{Sloane} } \\[2pt]\hline
(3,2) & \frac{C(t)}{\sqrt{1-4t} }            & A001700=\binom{2n+1}{n}                 & A111418 \  \mbox{\cite{Sloane} } \\[2pt]\hline
\end{array}
\end{eqnarray*}
Table 2. The specializations of $(x,y)$, where $C_{\frac{n}{2}}$ is
set to be zero when $n$ is odd.
\end{center}

\section{ Main results and bijective proofs }

\begin{lemma}\label{lemma 3.1.1}
There exists a bijection between the set $\mathscr{M}_{n, 0}(y+1,y)$
and the set $\bigcup_{\ell=0}^n\mathscr{M}_{n,\ell}(y,y)$.
\end{lemma}
\pf For any $P\in \mathscr{M}_{n, 0}(y+1,y)$, each $H$ step of $P$
on $X$-axis has weight $y+1$, or equivalently, it has weight $y$ or
$1$. If $P$ has $\ell$ $H$ steps weighted by $1$, replace each of
them by a $U$ step, we get a path $P^{*}\in
\mathscr{M}_{n,\ell}(y,y)$.

Conversely, for any $P^{*}\in \mathscr{M}_{n,\ell}(y,y)$, it has
exactly $\ell$ R-visible up steps, replace each of them by an $H$
step, we get a path $P$ has $\ell$ $H$ steps which are weighted by
$1$ and lying on $X$-axis.

Clearly, the above process indeed forms a bijection between the set
$\mathscr{M}_{n, 0}(y+1,y)$ and the set
$\bigcup_{\ell=0}^n\mathscr{M}_{n,\ell}(y,y)$. \qed\vskip0.2cm

\begin{theorem}
Let $\mathcal{M}=(M_{n,k}(x,y))_{n\geq k\geq 0}$ be given in Section
2. For any integers $n, r\geq 0$ and $m\geq \ell\geq 0$, set $N_r=\min\{n+r+1, m+r-\ell\}$.
Then there hold
\begin{eqnarray}
\sum_{k=0}^{N_r}\det\left(\begin{array}{cc}
M_{n,k}(x,y)   & M_{m,k+\ell+1}(x,y) \\[5pt]
M_{n+r+1,k}(x,y) & M_{m+r+1,k+\ell+1}(x,y)
\end{array}\right)
\hskip-.22cm &=&\hskip-.22cm \sum_{i=0}^{r} M_{n+i,0}(x,y)M_{m+r-i,\ell}(y,y), \label{eqn 3.1.1} \\
\sum_{\ell=0}^{m}\sum_{k=0}^{N_r}\det\left(\begin{array}{cc}
M_{n,k}(x,y)   & M_{m,k+\ell+1}(x,y) \\[5pt]
M_{n+r+1,k}(x,y) & M_{m+r+1,k+\ell+1}(x,y)
\end{array}\right)
\hskip-.22cm &=&\hskip-.22cm \sum_{i=0}^{r} M_{n+i,0}(x,y)M_{m+r-i,0}(y+1, y).\label{eqn 3.1.2}
\end{eqnarray}
\end{theorem}

\pf Define
\begin{eqnarray*}
\mathscr{A}_{n,m,k,\ell}^{(r)}(x,y) \hskip-.22cm &=&\hskip-.22cm \{(P, Q)|P\in \mathscr{M}_{n, k}(x,y), Q\in \mathscr{M}_{m+r+1, k+\ell+1}(x,y)\}, \\
\mathscr{B}_{n,m,k,\ell}^{(r)}(x,y) \hskip-.22cm &=&\hskip-.22cm \{(P, Q)|P\in \mathscr{M}_{n+r+1, k}(x,y), Q\in \mathscr{M}_{m, k+\ell+1}(x,y)\},
\end{eqnarray*}
and $\mathscr{C}_{n,m,k,\ell}^{(r,i)}(x,y)$ to be the subset of $\mathscr{A}_{n,m,k,\ell}^{(r)}(x,y)$ such that for any
$(P, Q)\in \mathscr{C}_{n,m,k,\ell}^{(r,i)}(x,y)$ $Q=Q_{1}UQ_{2}$ with $Q_{1}\in \mathscr{M}_{i,k}(x,y)$ and $Q_2\in
\mathscr{M}_{m+r-i,\ell}(y,y)$ for $k\leq i\leq r$. In other words, for any
$(P, Q)\in \mathscr{C}_{n,m,k,\ell}^{(r,i)}(x,y)$, $Q$ satisfies the conditions that
(a) the last $(\ell+1)$-th R-visible up step of $Q$ stays at level $k+1$, and (b) there are
exactly $i$ steps immediately ahead of the last $(\ell+1)$-th R-visible up step of
$Q$ for $k\leq i\leq r$.

It is clear that the weights of the sets
$\mathscr{A}_{n,m,k,\ell}^{(r)}(x,y)$ and $\mathscr{B}_{n,m,k,\ell}^{(r)}(x,y)$ are
\begin{eqnarray*}
w(\mathscr{A}_{n,m,k,\ell}^{(r)}(x,y)) \hskip-.22cm &=&\hskip-.22cm M_{n, k}(x,y)M_{m+r+1, k+\ell+1}(x,y), \\
w(\mathscr{B}_{n,m,k,\ell}^{(r)}(x,y)) \hskip-.22cm &=&\hskip-.22cm M_{n+r+1, k}(x,y)M_{m, k+\ell+1}(x,y).
\end{eqnarray*}
Given $0\leq i\leq r$, the weight of the set $\bigcup_{k=0}^{i}\mathscr{C}_{n,m,k,\ell}^{(r,i)}(x,y)$
is $M_{n+i,0}(x,y)M_{m+r-i,\ell}(y,y)$. In fact, this claim can be verified by the following argument.
For any $(P, Q)\in \mathscr{C}_{n,m,k,\ell}^{(r,i)}(x,y)$,
we have $Q=Q_{1}UQ_{2}$ as mentioned above with $Q_{1}\in \mathscr{M}_{i,k}(x,y)$ and $Q_2\in
\mathscr{M}_{m+r-i,\ell}(y,y)$, then $P\overline{Q}_{1}\in \mathscr{M}_{n+i,0}(x,y)$ such that the last
$(i+1)$-th step of $P\overline{Q}_{1}$ is at level $k$. Summing $k$ for $0\leq k\leq i$, all
$P\overline{Q}_{1}\in \mathscr{M}_{n+i,0}(x,y)$ contribute the total weight $M_{n+i,0}(x,y)$ and all $Q_2\in
\mathscr{M}_{m+r-i,\ell}(y,y)$ contribute the total weight $M_{m+r-i,\ell}(y,y)$. Hence,
$w(\bigcup_{k=0}^{i}\mathscr{C}_{n,m,k,\ell}^{(r,i)}(x,y))=M_{n+i,0}(x,y)M_{m+r-i,\ell}(y,y)$,
and then
\begin{eqnarray*}
w(\bigcup_{i=0}^{r}\bigcup_{k=0}^{i}\mathscr{C}_{n,m,k,\ell}^{(r,i)}(x,y))
=w(\bigcup_{k=0}^{r}\bigcup_{i=k}^{r}\mathscr{C}_{n,m,k,\ell}^{(r,i)}(x,y))=\sum_{i=0}^{r} M_{n+i,0}(x,y)M_{m+r-i,\ell}(y,y).
\end{eqnarray*}

Let $\mathscr{A}_{n,m,\ell}^{(r)}(x,y)=\bigcup_{k=0}^{N_r}\mathscr{A}_{n,m,k,\ell}^{(r)}(x,y)$,
$\mathscr{B}_{n,m,\ell}^{(r)}(x,y)=\bigcup_{k=0}^{N_r}\mathscr{B}_{n,m,k,\ell}^{(r)}(x,y)$
and $\mathscr{C}_{n,m,\ell}^{(r)}(x,y)=\bigcup_{0\leq k\leq i\leq r}\mathscr{C}_{n,m,k,\ell}^{(r,i)}(x,y)$.
In order to prove (\ref{eqn 3.1.1}), it suffices to construct a simple bijection $\phi$ between
$\mathscr{A}_{n,m,\ell}^{(r)}(x,y)-\mathscr{C}_{n,m,\ell}^{(r)}(x,y)$ and
$\mathscr{B}_{n,m,\ell}^{(r)}(x,y)$ such that the $\phi$ is still preserving the
weights.

For any $(P, Q)\in
\mathscr{A}_{n,m,k,\ell}^{(r)}(x,y)-\bigcup_{i=k}^{r}\mathscr{C}_{n,m,k,\ell}^{(r,i)}(x,y)$, $Q$ has exactly
$k+\ell+1$ R-visible up steps. We claim that there always exists a path $Q'$ of length $r+1$
which is immediately ahead of the last $(\ell+1)$-th (also the $(k+1)$-th along the path) R-visible up step
of $Q$. Otherwise, $Q'\in \mathscr{M}_{i,k}(x,y)$ for some $k\leq i\leq r$, and then
$(P, Q)\in \bigcup_{i=k}^{r}\mathscr{C}_{n,m,k,\ell}^{(r,i)}(x,y)$, a contradiction.

For any $(P, Q)\in \mathscr{A}_{n,m,\ell}^{(r)}(x,y)-\mathscr{C}_{n,m,\ell}^{(r)}(x,y)$,
find the path $Q'$ of length $r+1$ which is immediately ahead of the last $(\ell+1)$-th
R-visible up step of $Q$, namely, $Q$ can be uniquely partitioned into $Q=Q_1Q'U^{*}Q_2$,
where $Q_1\in \mathscr{M}_{j}(x,y)$ for some $j\geq 0$, $Q_2\in \mathscr{M}_{\ell}(y,y)$ and
$U^{*}$ is the last $(\ell+1)$-th R-visible up step of $Q$. Then
we can construct $\phi(P, Q)=(P^*,
Q^{*})\in \mathscr{B}_{n,m,\ell}^{(r)}(x,y)$ as follows: (1) delete
the path $Q'$ in $Q$ to get $Q^*=Q_1U^{*}Q_2$; (2) annex the reverse path $\overline{Q'}$ of $Q'$ to the end of $P$ to get
$P^{*}$, that is, $P^*=P\overline{Q'}$. More precisely,

\begin{itemize}
\item{} $(P, Q)\in \mathscr{A}_{n,m,k,\ell}^{(r)}(x,y)-\bigcup_{i=k}^{r}\mathscr{C}_{n,m,k,\ell}^{(r,i)}(x,y)$
leads to $(P^*, Q^{*})\in \mathscr{B}_{n,m,j,\ell}^{(r)}(x,y)$, where $Q=Q_1Q'U^{*}Q_2$ is factored as above.
\end{itemize}

Note that in this case the last $(\ell+1)$-th R-visible up step of $Q$ are still the
one of $Q^{*}$.

Conversely, we can recover $(P, Q)\in
\mathscr{A}_{n,m,\ell}^{(r)}(x,y)-\mathscr{C}_{n,m,\ell}^{(r)}(x,y)$
from $(P^*, Q^{*})\in \mathscr{B}_{n,m,\ell}^{(r)}(x,y)$.
For any $(P^*, Q^{*})\in \mathscr{B}_{n,m,\ell}^{(r)}(x,y)$, $P^{*}$ can be uniquely
partitioned into $P^{*}=PP'$ such that $P\in \mathscr{M}_{n, k}(x,y)$ for some $0\leq k\leq N_r$
and $P'$ has length $r+1$. Then delete the path $P'$ of $P^*$ to get $P$,
and interpolate the reverse path $\overline{P'}$ of $P'$ immediately ahead of the
last $(\ell+1)$-th R-visible up step of $Q^{*}$ to get $Q$. In this case, the last $(\ell+1)$-th R-visible
up step of $Q^{*}$ are also the one of $Q$ and there are at least $r+1$ steps immediately ahead of it.

Note that $\phi$ does not change the weight of any $H$ step,
despite $\phi$ possibly exchange some $U$ steps and $D$ steps,
but all $U$ and $D$ steps have the same weight $1$. Hence,
$\phi$ is indeed a bijection and also preserves weights.
Therefore, (\ref{eqn 3.1.1}) is proved.

Summing the two sides of (\ref{eqn 3.1.1}) for $0\leq \ell\leq m+r$, by
Lemma 3.1, (\ref{eqn 3.1.2}) follows. \qed\vskip0.2cm

The special case $r=0$ in (\ref{eqn 3.1.1}) produces the following result.
\begin{theorem}
Let $\mathcal{M}=(M_{n,k}(x,y))_{n\geq k\geq 0}$ be given in Section
2. For any integers $n\geq 0$ and $m\geq \ell\geq 0$, set $N_0=\min\{n+1, m-\ell\}$.
Then there hold
\begin{eqnarray}\label{eqn 3.2.1}
\sum_{k=0}^{N_0}\det\left(\begin{array}{cc}
M_{n,k}(x,y)   & M_{m,k+\ell+1}(x,y) \\[5pt]
M_{n+1,k}(x,y) & M_{m+1,k+\ell+1}(x,y)
\end{array}\right)
\hskip-.22cm &=&\hskip-.22cm M_{n,0}(x,y)M_{m,\ell}(y,y).
\end{eqnarray}
\end{theorem}

Now we concentrate on the specialization of the parameters $(x,y)$ in Theorem 3.3, which generates
many identities involving Catalan numbers.

\subsubsection*{\bf Example (i)} When $(x, y)=(1,2)$, (\ref{eqn 2.1}) and (\ref{eqn 2.2}) yield that
$M(2,2;t)=C^2(t)$ and $M(1,2;t)=C(t)$, so
$\mathcal{M}=(M_{n,k}(1,2))_{n\geq k\geq 0}$ is the Riordan array
$(C(t), tC^2(t))$. By the series expansion \cite{Stanley},
\begin{eqnarray}\label{eqn 3.3.1}
C(t)^{\alpha} \hskip-.22cm &=&\hskip-.22cm \sum_{n\geq
0}\frac{\alpha}{2n+\alpha}\binom{2n+\alpha}{n} t^n,
\end{eqnarray}
we have
\begin{eqnarray}\label{eqn 3.3.2}
M_{n,k}(1,2) \hskip-.22cm &=&\hskip-.22cm
[t^n]C(t)(tC^2(t))^{k}=[t^{n-k}]C(t)^{2k+1}=\frac{2k+1}{2n+1}\binom{2n+1}{n-k}.
\end{eqnarray}
Then, after some routine simplifications, (\ref{eqn 3.2.1}) produces the following result.
\begin{eqnarray}\label{eqn 3.3.3}
\frac{\ell+1}{m+1}\binom{2m+2}{m-\ell}C_{n}\hskip-.22cm
&=&\hskip-.22cm
\sum_{k=0}^{N_0}\frac{(2k+1)(2k+2\ell+1)\alpha_{n,k}(m,\ell)}{(2n+1)_3(2m+1)_3}\binom{2n+3}{n-k+1}\binom{2m+3}{m-k-\ell}, \hskip-1cm
\end{eqnarray}
where $\alpha_{n,k}(m,\ell)=6(m-n)(n+1)(m+1)+(\ell+1)(2k+\ell+2)(2n+1)(2n+2)-2(m-n)k(k+1)(2n+2m+3)$ and  $(x)_k=x(x+1)\cdots(x+k-1)$ for $k\geq 1$
and $(x)_0=1$.

Taking $\ell=0$ and $m=n-1, n$ or $n+1$ into account, we have
\begin{eqnarray*}
\alpha_{n,k}(n-1,0) \hskip-.22cm &=&\hskip-.22cm (n+k+3)(8nk+2n+2k+2),\\
\alpha_{n,k}(n,0)   \hskip-.22cm &=&\hskip-.22cm (2k+2)(2n+1)(2n+2),\\
\alpha_{n,k}(n+1,0) \hskip-.22cm &=&\hskip-.22cm (n-k+1)(8nk+14n+10k+16).
\end{eqnarray*}
Then in these three cases, after shifting $n$ to $n+1$ in the case $m=n-1$, (\ref{eqn 3.3.3}) generates
\begin{corollary} For any integer $n\geq 0$, there hold
\begin{eqnarray}
C_{n+1}^2\hskip-.22cm &=&\hskip-.22cm
\sum_{k=0}^{n}\frac{(2k+1)(2k+3)(8nk+2n+10k+4)}{(2n+1)(2n+2)(2n+3)(2n+4)(2n+5)}\binom{2n+2}{n-k}\binom{2n+5}{n-k+2}, \nonumber\\
C_{n}C_{n+1} \hskip-.22cm &=&\hskip-.22cm
\sum_{k=0}^{n}\frac{(2k+1)(2k+2)(2k+3)}{(2n+1)(2n+2)(2n+3)^2}\binom{2n+3}{n-k}\binom{2n+3}{n-k+1}, \label{eqn 3.3.4}\\
C_{n}C_{n+2}\hskip-.22cm &=&\hskip-.22cm
\sum_{k=0}^{n}\frac{(2k+1)(2k+3)(8nk+14n+10k+16)}{(2n+1)(2n+2)(2n+3)(2n+4)(2n+5)}\binom{2n+2}{n-k}\binom{2n+5}{n-k+1}.
\nonumber
\end{eqnarray}
\end{corollary}
\vskip0.2cm

\subsubsection*{\bf Example (ii)}  When $(x, y)=(2,2)$, (\ref{eqn 2.2}) yields that
$M(2,2;t)=C^2(t)$, so $\mathcal{M}=(M_{n,k}(2,2))_{n\geq k\geq 0}$
is the Riordan array $(C^2(t), tC^2(t))$, it is also Shapiro's Catalan triangle aforementioned. By (\ref{eqn 3.3.2}), we have
\begin{eqnarray}\label{eqn 3.4.1}
M_{n,k}(2,2) \hskip-.22cm &=&\hskip-.22cm
[t^n]C^2(t)(tC^2(t))^{k}=[t^{n-k}]C(t)^{2k+2}=\frac{2k+2}{2n+2}\binom{2n+2}{n-k}.
\end{eqnarray}
Then, after some routine simplifications, (\ref{eqn 3.2.1}) produces the following result.
\begin{eqnarray}\label{eqn 3.4.2}
\frac{\ell+1}{m+1}\binom{2m+2}{m-\ell}C_{n+1} \hskip-.22cm
&=&\hskip-.22cm
\sum_{k=0}^{N_0}\frac{(2k+2)(2k+2\ell+4)\beta_{n,k}(m,\ell)}{(2n+2)_3(2m+2)_3}\binom{2n+4}{n-k+1}\binom{2m+4}{m-k-\ell}, \hskip-1cm
\end{eqnarray}
where $\beta_{n,k}(m,\ell)=6(m-n)(n+1)(m+1)+(\ell+1)(2k+\ell+3)(2n+2)(2n+3)-2(m-n)k(k+2)(2n+2m+5)$.

Taking $\ell=0$ and $m=n-1, n$ or $n+1$ into account, we have
\begin{eqnarray*}
\beta_{n,k}(n-1,0) \hskip-.22cm &=&\hskip-.22cm (n+k+3)(8nk+6n+6k+6), \\
\beta_{n,k}(n,0)   \hskip-.22cm &=&\hskip-.22cm (2k+3)(2n+2)(2n+3),\\
\beta_{n,k}(n+1,0) \hskip-.22cm &=&\hskip-.22cm (n-k+1)(8nk+18n+14k+30).
\end{eqnarray*}
Then in these three cases, after shifting $n$ to $n+1$ in the case $m=n-1$, (\ref{eqn 3.4.2}) generates
\begin{corollary} For any integer $n\geq 0$, there hold
\begin{eqnarray}
C_{n+1}C_{n+2}\hskip-.22cm &=&\hskip-.22cm
\sum_{k=0}^{n}\frac{(2k+2)(2k+4)(8nk+6n+14k+12)}{(2n+2)(2n+3)(2n+4)(2n+5)(2n+6)}\binom{2n+3}{n-k}\binom{2n+6}{n-k+2}, \nonumber\\
C_{n+1}^2 \hskip-.22cm &=&\hskip-.22cm
\sum_{k=0}^{n}\frac{(2k+2)(2k+3)(2k+4)}{(2n+2)(2n+3)(2n+4)^2}\binom{2n+4}{n-k}\binom{2n+4}{n-k+1}, \label{eqn 3.4.3}\\
C_{n+1}C_{n+2}\hskip-.22cm &=&\hskip-.22cm
\sum_{k=0}^{n}\frac{(2k+2)(2k+4)(8nk+18n+14k+30)}{(2n+2)(2n+3)(2n+4)(2n+5)(2n+6)}\binom{2n+3}{n-k}\binom{2n+6}{n-k+1}.
\nonumber
\end{eqnarray}
\end{corollary}

\subsubsection*{\bf Example (iii)}  When $(x, y)=(3,2)$, (\ref{eqn 2.1}) and (\ref{eqn 2.2}) yield that
$M(2,2;t)=C^2(t)$ and $M(3,2;t)=\frac{C(t)}{\sqrt{1-4t}}$, so
$\mathcal{M}=(M_{n,k}(3,2))_{n\geq k\geq 0}$ is the Riordan array
$(\frac{C(t)}{\sqrt{1-4t}}, tC^2(t))$. By the series expansion
\cite{Stanley},
\begin{eqnarray*}
\frac{C(t)^{\alpha}}{\sqrt{1-4t}} \hskip-.22cm &=&\hskip-.22cm
\sum_{n\geq 0}\binom{2n+\alpha}{n} t^n,
\end{eqnarray*}
we have
\begin{eqnarray}\label{eqn 3.5.1}
M_{n,k}(3,2) \hskip-.22cm &=&\hskip-.22cm
[t^n]\frac{C(t)}{\sqrt{1-4t}}(tC^2(t))^{k}=[t^{n-k}]\frac{C(t)^{2k+1}}{\sqrt{1-4t}}=\binom{2n+1}{n-k}.
\end{eqnarray}
Then, after some routine simplifications, (\ref{eqn 3.2.1}) produces the following result.
\begin{eqnarray}\label{eqn 3.5.2}
\frac{\ell+1}{m+1}\binom{2m+2}{m-\ell}\binom{2n+1}{n} \hskip-.22cm
&=&\hskip-.22cm
\sum_{k=0}^{N_0}\frac{\gamma_{n,k}(m,\ell)}{(2n+2)_2(2m+2)_2}\binom{2n+3}{n-k+1}\binom{2m+3}{m-k-\ell}, \hskip-.5cm
\end{eqnarray}
where $\gamma_{n,k}(m)=2(m-n)(n+1)(m+1)+(\ell+1)(2k+\ell+2)(2n+2)(2n+3)-2(m-n)k(k+1)(2n+2m+5)$.

Taking $\ell=0$ and $m=n-1, n$ or $n+1$ into account, we have
\begin{eqnarray*}
\gamma_{n,k}(n-1,0)\hskip-.22cm &=&\hskip-.22cm (n+k+2)(8nk+6n+6k+6), \\
\gamma_{n,k}(n,0)  \hskip-.22cm &=&\hskip-.22cm (2k+2)(2n+2)(2n+3), \\
\gamma_{n,k}(n+1,0)\hskip-.22cm &=&\hskip-.22cm (n-k+1)(8nk+10n+14k+6).
\end{eqnarray*}
Then in these three cases, after shifting $n$ to $n+1$ in the case $m=n-1$, (\ref{eqn 3.4.2}) generates
\begin{corollary}
For any integer $n\geq 0$, there hold
\begin{eqnarray}
\binom{2n+3}{n+1}C_{n+1} \hskip-.22cm &=&\hskip-.22cm
\sum_{k=0}^{n}\frac{(8nk+6n+14k+12)}{(2n+2)(2n+3)(2n+4)}\binom{2n+3}{n-k}\binom{2n+4}{n-k+2}, \nonumber \\
\binom{2n+1}{n}C_{n+1} \hskip-.22cm &=&\hskip-.22cm
\sum_{k=0}^{n}\frac{(2k+2)}{(2n+2)(2n+3)}\binom{2n+3}{n-k}\binom{2n+3}{n-k+1}, \label{eqn 3.5.3} \\
\binom{2n+1}{n}C_{n+2} \hskip-.22cm &=&\hskip-.22cm
\sum_{k=0}^{n}\frac{(8nk+10n+14k+6)}{(2n+2)(2n+3)(2n+4)}\binom{2n+4}{n-k}\binom{2n+3}{n-k+1}.
\nonumber
\end{eqnarray}
\end{corollary}

\subsubsection*{\bf Example (iv)} When $(x, y)=(0,0)$, (\ref{eqn 2.2}) yields that $M(0,0;t)=C(t^2)$,
so $\mathcal{M}=(M_{n,k}(0,0))_{n\geq k\geq 0}$ is the Riordan array
$(C(t^2), tC(t^2))$. By (\ref{eqn 3.3.2}), we have
\begin{eqnarray}\label{eqn 3.6.1}
\hskip1cm M_{n,k}(0,0) \hskip-.22cm &=&\hskip-.22cm
[t^n]C(t^2)(tC(t^2))^{k}=[t^{n-k}] C(t^2)^{k+1}=\left\{
\begin{array}{ll}
\frac{k+1}{n+1}\binom{n+1}{\frac{n-k}{2}}, & \mbox{if}\ n-k\ \mbox{ even }, \\
0, & \mbox{otherwise}.
\end{array}\right.
\end{eqnarray}
Replacing $n, m, \ell$ by $2n, 2m, 2\ell$ respectively, after some routine simplifications, (\ref{eqn 3.2.1}) produces the following result.
\begin{eqnarray}\label{eqn 3.6.2}
\frac{2\ell+1}{2m+1}\binom{2m+1}{m-\ell}C_n \hskip-.22cm
&=&\hskip-.22cm
\sum_{k=0}^{N_0}\frac{\lambda_{n,k}(m,\ell)}{(2n+1)_2(2m+1)_2}\binom{2n+2}{n-k}\binom{2m+2}{m-k-\ell},
\end{eqnarray}
where
$\lambda_{n,k}(m,\ell)=2(2n+1)(2k+\ell+2)(2(k+1)(k+\ell+1)-(m+1))+2(m-n)(k+\ell+1)(2k+1)(2k+3)$.


Taking $\ell=0$ and $n=m$ into account, we have
\begin{eqnarray*}
\lambda_{m,k}(m,0)\hskip-.22cm &=&\hskip-.22cm (2m+1)(2k+2)((2k+1)(2k+3)-(2m+1)).
\end{eqnarray*}
In this case, (\ref{eqn 3.6.2}) generates
\begin{corollary} For any integer $m\geq 0$, there holds
\begin{eqnarray}\label{eqn 3.6.3}
C_{m}^2 \hskip-.22cm &=&\hskip-.22cm
\sum_{k=0}^{m}\frac{(2k+2)((2k+1)(2k+3)-(2m+1))}{(2m+1)(2m+2)^2}\binom{2m+2}{m-k}^2.
\end{eqnarray}
\end{corollary}

\begin{remark}
It should be pointed out that despite (\ref{eqn 3.2.1}) is not valid
for any integer $\ell\leq -1$, since in these cases $M_{m, \ell}(y,
y)$ has not be defined, but (\ref{eqn 3.3.3}), (\ref{eqn 3.4.2}), (\ref{eqn 3.5.2}), (\ref{eqn 3.6.2}) are
all correct for any integer $\ell\leq -1$ if one notices that they
hold trivially for any integer $\ell> m$ and both sides of them can
be transferred into polynomials on $\ell$.
\end{remark}

It is also worth pointing out that the case $m=n$ in (\ref{eqn 3.3.3}), (\ref{eqn 3.4.2}) and (\ref{eqn 3.5.2}) generates the following
results respectively.

\begin{corollary} For any integers $n\geq\ell\geq 0$, there hold
\begin{eqnarray}
\frac{1}{n+1}\binom{2n+2}{n-\ell}C_n \hskip-.22cm &=&\hskip-.22cm
\sum_{k=0}^{n-\ell}
\frac{(2k+1)(2k+\ell+2)(2k+2\ell+3)}{(2n+1)(2n+2)(2n+3)^2}\binom{2n+3}{n-k-\ell}\binom{2n+3}{n-k+1}, \hskip-2cm \label{eqn 3.7.1}\\
\frac{1}{n+1}\binom{2n+2}{n-\ell}C_{n+1} \hskip-.22cm
&=&\hskip-.22cm
\sum_{k=0}^{n-\ell}\frac{(2k+2)(2k+\ell+3)(2k+2\ell+4)}{(2n+2)(2n+3)(2n+4)^2}\binom{2n+4}{n-k-\ell}\binom{2n+4}{n-k+1},\hskip-2cm  \label{eqn 3.7.2}\\
\frac{1}{n+1}\binom{2n+2}{n-\ell}\binom{2n+1}{n} \hskip-.22cm
&=&\hskip-.22cm
\sum_{k=0}^{n-\ell}\frac{(2k+\ell+2)}{(2n+2)(2n+3)}\binom{2n+3}{n-k-\ell}\binom{2n+3}{n-k+1}.
\hskip-2cm  \label{eqn 3.7.3}
\end{eqnarray}
\end{corollary}

\begin{remark}
Recently, Guti¨¦rrez et al. \cite{Gutierrez}, Miana and Romero \cite{Miana},
Chen and Chu \cite{ChenChu}, and Guo and Zeng \cite{GuoZeng} studied the binomial
sums $\sum_{k=0}^{n}(k+1)^m\binom{2n+2}{n-k}^2$ related to the classical
Catalan triangle \cite{ShapA}. Zhang and Pang \cite{ZhangPang} also considered some alternating cases.
Later, Miana and Romero \cite{MianaRom} investigated another binomial
sums $\sum_{k=0}^{n}(2k+1)^m\binom{2n+1}{n-k}^2$. Clearly our work is closely related to theirs
from a different direction. Setting $p=k+1, \ell=n-i+1$, and then replacing $n$ by $n-2$, (\ref{eqn 3.7.2}) reduces
to the main identity obtained by Guti¨¦rrez et al. \cite[Theorem 5]{Gutierrez}.
As mentioned in Remark 3.7, (\ref{eqn 3.7.1})-(\ref{eqn 3.7.3}) also hold for any integer $\ell<0$. Specially, in the case
$\ell=-1$, replacing $n+1$ by $n$, after some routine
simplifications, (\ref{eqn 3.7.1})-(\ref{eqn 3.7.3}) lead respectively
to the following identities,
\begin{eqnarray}
\sum_{k=0}^{n}\frac{(2k+1)^3}{(2n+1)^2}\binom{2n+1}{n-k}^2\hskip-.22cm &=&\hskip-.22cm
\binom{2n}{n}^2,\hskip1cm \mbox{\cite[Remark 11]{MianaRom}} \label{eqn 3.9E}\\
\sum_{k=0}^{n}\frac{(k+1)^3}{(n+1)^2}\binom{2n+2}{n-k}^2 \hskip-.22cm &=&\hskip-.22cm
\binom{2n}{n}\binom{2n+1}{n}, \hskip1cm\mbox{\cite[Corollary 6]{Gutierrez}} \label{eqn 3.9F}\\
\sum_{k=0}^{n}\frac{(2k+1)}{(2n+1)}\binom{2n+1}{n-k}^2 \hskip-.22cm
&=&\hskip-.22cm \binom{2n}{n}^2. \label{eqn 3.9G}
\end{eqnarray}
Note that (\ref{eqn 3.9E}), (\ref{eqn 3.9G}) and (\ref{eqn 4.2B})
can be regarded as companion ones of an identity obtained by Deng
and Yan \cite{Deng},
\begin{eqnarray*}
\sum_{k=0}^{n}\frac{(2k+1)^2}{(2n+1)}\binom{2n+1}{n-k} \hskip-.22cm
&=&\hskip-.22cm 4^n.
\end{eqnarray*}
Moreover, (\ref{eqn 3.7.1})-(\ref{eqn 3.7.3}) hold symmetrically on
$\ell=-1$. For example, in the case $\ell=-2$, (\ref{eqn 3.7.1})-(\ref{eqn 3.7.3})
can reduce respectively to (\ref{eqn 3.3.4})
(\ref{eqn 3.4.3}) and (\ref{eqn 3.5.3}) which also correspond to the
case $\ell=0$ in (\ref{eqn 3.7.1})-(\ref{eqn 3.7.3}).
\end{remark}

The special case $r=1$, $\ell=0$ and $m=n$ in (\ref{eqn 3.1.1}) produces the following result.
\begin{theorem}
Let $\mathcal{M}=(M_{n,k}(x,y))_{n\geq k\geq 0}$ be given in Section
2. Then there holds
\begin{eqnarray}\label{eqn 3.10}
\sum_{k=0}^{n}\det\left(\begin{array}{cc}
M_{n,k}(y,y) & M_{n,k+1}(y,y) \\[5pt]
M_{n+2,k}(y,y) & M_{n+2,k+1}(y,y)
\end{array}\right)
\hskip-.22cm &=&\hskip-.22cm 2M_{n,0}(y,y)M_{n+1,0}(y,y),
\end{eqnarray}
\end{theorem}

In the case $y=2$, together with (\ref{eqn 3.4.1}), after some routine
computations, (\ref{eqn 3.10}) generates

\begin{corollary} For any integer $n\geq 0$, there holds
\begin{eqnarray*}
C_nC_{n+1} \hskip-.22cm &=&\hskip-.22cm \sum_{k=0}^{n}
\frac{(2k+2)(2k+3)(2k+4)}{(2n+2)(2n+3)(2n+6)(2n+7)}\binom{2n+3}{n-k}\binom{2n+7}{n-k+2}.
\end{eqnarray*}
\end{corollary}

\vskip1cm
\section{ Alternating Cases}

In this section, we consider some alternating sums related to
$\mathcal{M}=(M_{n,k}(x,y))_{n\geq k\geq 0}$. Despite it has no
general and unified results as in the previous section, but in
several isolated cases, mainly by the creative telescoping algorithm
\cite{Petkov, Zeilberger}, we also obtain some interesting results.

\begin{theorem}
Let $\mathcal{M}=(M_{n,k}(x,y))_{n\geq k\geq 0}$ be given in Section
2. Then there holds
\begin{eqnarray}\label{eqn 4.1}
\sum_{k=0}^{n}(-1)^{n-k}\det\left(\begin{array}{cc}
M_{n,k}(0,0) & M_{n,k+1}(0,0) \\[5pt]
M_{n+1,k}(0,0) & M_{n+1,k+1}(0,0)
\end{array}\right)
\hskip-.22cm &=&\hskip-.22cm C_{n+1},
\end{eqnarray}
or equivalently,
\begin{eqnarray}
C_{2m+1} \hskip-.22cm &=&\hskip-.22cm
\sum_{j=0}^{m}\frac{(2j+2)^2}{(2m+2)^2}\binom{2m+2}{m-j}^2, \label{eqn 4.2}\\
C_{2m+2} \hskip-.22cm &=&\hskip-.22cm
\sum_{j=0}^{m}\frac{(2j+2)^2}{(2m+2)(2m+3)}\binom{2m+3}{m-j}\binom{2m+3}{m-j+1}.
\label{eqn 4.3}
\end{eqnarray}
\end{theorem}

\pf By (\ref{eqn 3.6.1}), the cases $n=2m$ and $n=2m+1$ in (\ref{eqn
4.1}) are equivalent to (\ref{eqn 4.2}) and (\ref{eqn 4.3})
respectively.

For (\ref{eqn 4.2}), its right side counts the set of pairs $(P,Q)$
such that $P, Q\in \mathscr{M}_{m, j}(2,2)$ for $0\leq j\leq m$. Let
$\overline{Q}$ be the reverse path of $Q$.
Clearly, $P\overline{Q}\in \mathscr{M}_{2m, 0}(2,2)$. The process is
obviously reversible. This builds a simple bijection between the
sets $\bigcup_{0\leq j\leq m}\big(\mathscr{M}_{m, j}(2,2)\times
\mathscr{M}_{m, j}(2,2)\big)$ and $\mathscr{M}_{2m, 0}(2,2)$, which,
by (\ref{eqn 3.4.1}), proves (\ref{eqn 4.2}).

Note that $(j+1)^2-(m+2)^2=-(m+j+3)(m-j+1)$, by the trivial
identities
\begin{eqnarray*}
\sum_{j=0}^{m}\binom{2m+3}{m-j}\binom{2m+3}{m-j+1} \hskip-.22cm &=&\hskip-.22cm \sum_{j=-m-2}^{-2}\binom{2m+3}{m-j}\binom{2m+3}{m-j+1}, \\
\sum_{j=0}^{m}\binom{2m+2}{m-j}^2 \hskip-.22cm &=&\hskip-.22cm
\sum_{j=-m-2}^{-2}\binom{2m+2}{m-j}^2,
\end{eqnarray*}
we have
\begin{eqnarray*}
\lefteqn{\sum_{j=0}^{m}\frac{(2j+2)^2}{(2m+2)(2m+3)}\binom{2m+3}{m-j}\binom{2m+3}{m-j+1} } \\
\hskip-.22cm &=&\hskip-.22cm 4 \sum_{j=0}^{m}\frac{(m+2)^2}{(2m+2)(2m+3)}\binom{2m+3}{m-j}\binom{2m+3}{m-j+1} \\
& & \hskip1cm +4\sum_{j=0}^{m}\frac{(j+1)^2-(m+2)^2}{(2m+2)(2m+3)}\binom{2m+3}{m-j}\binom{2m+3}{m-j+1} \\
\hskip-.22cm &=&\hskip-.22cm
\frac{4(m+2)^2}{(2m+2)(2m+3)}\sum_{j=0}^{m}\binom{2m+3}{m-j}\binom{2m+3}{m-j+1}
-\frac{4(2m+3)}{(2m+2)}\sum_{j=0}^{m}\binom{2m+2}{m-j}^2  \\
\hskip-.22cm &=&\hskip-.22cm \frac{2(m+2)^2}{(2m+2)(2m+3)}\left\{\sum_{j=-m-2}^{m}\binom{2m+3}{m-j}\binom{2m+3}{m-j+1}-\binom{2m+3}{m+1}^2\right\} \\
& & \hskip1cm -\frac{2(2m+3)}{(2m+2)}\left\{\sum_{j=-m-2}^{m}\binom{2m+2}{m-j}^2-\binom{2m+2}{m+1}^2 \right\} \\
\hskip-.22cm &=&\hskip-.22cm \frac{2(m+2)^2}{(2m+2)(2m+3)}\left\{\binom{4m+6}{2m+2}-\binom{2m+3}{m+1}^2\right\} \\
& & \hskip1cm -\frac{2(2m+3)}{(2m+2)}\left\{\binom{4m+4}{2m+2}-\binom{2m+2}{m+1}^2 \right\} \ \ (\mbox{by Vandermonde identity}) \\
\hskip-.22cm &=&\hskip-.22cm C_{2m+2},
\end{eqnarray*}
which proves (\ref{eqn 4.3}).\qed

\begin{remark}
Note that $\mathcal{M}=(M_{m,k}(x,y))_{m\geq k\geq 0}$ also forms an
admissible matrix defined by Aigner \cite{Aigner}, then (\ref{eqn
4.2}) can be obtained easily by the elementary property of the
admissible matrix $(M_{m,k}(2,2))_{m\geq k\geq 0}$. What's more, one
can also get another one from the admissible matrix
$(M_{m,k}(1,2))_{m\geq k\geq 0}$,
\begin{eqnarray}\label{eqn 4.2B}
\sum_{j=0}^{m}\frac{(2j+1)^2}{(2m+1)^2}\binom{2m+1}{m-j}^2
\hskip-.22cm &=&\hskip-.22cm  C_{2m},
\end{eqnarray}
which also has a combinatorial proof similar to (\ref{eqn 4.2}).
Using the creative telescoping algorithm, one can find another identity
related to (\ref{eqn 3.6.3}) and (\ref{eqn 4.2}), that is,
\begin{eqnarray}\label{eqn 4.2C}
\sum_{j=0}^{m}\frac{(j+1)}{(m+1)}\binom{2m+2}{m-j}^2 \hskip-.22cm
&=&\hskip-.22cm \binom{2m+1}{m}^2,
\end{eqnarray}
which is also obtained implicitly by Chen and Chu \cite[Corollary 4]{ChenChu}.

Note that (\ref{eqn 3.9F}), (\ref{eqn 4.2}) and (\ref{eqn 4.2C}) can
also be regarded as companion ones of an identity obtained by
Cameron and Nkwanta \cite{CamNk},
\begin{eqnarray*}
\sum_{j=0}^{m}\frac{(j+1)^2}{(m+1)}\binom{2m+2}{m-j} \hskip-.22cm
&=&\hskip-.22cm 4^m.
\end{eqnarray*}
\end{remark}

\begin{theorem} For any integer $m\geq 0$, there hold
\begin{eqnarray}
\binom{2m+1}{m}   \hskip-.22cm &=&\hskip-.22cm \sum_{j=0}^{m}(-1)^{j}\frac{(2j+2)^2}{(2m+2)^2}\binom{2m+2}{m-j}^2,
\hskip1cm\mbox{\cite[Theorem 2.2]{ZhangPang}} \label{eqn 4.4}\\
\binom{2m+2}{m+1} \hskip-.22cm &=&\hskip-.22cm \sum_{j=0}^{m}(-1)^{j}\frac{(2j+2)^2}{(2m+2)(2m+3)}\binom{2m+3}{m-j}\binom{2m+3}{m-j+1}, \label{eqn 4.5} \\
\binom{2m}{m} \hskip-.22cm &=&\hskip-.22cm
\sum_{j=0}^{m}(-1)^{j}\frac{(2j+1)}{(2m+1)}\binom{2m+1}{m-j}^2.\label{eqn
4.6}
\end{eqnarray}
\end{theorem}
\pf For $0\leq j\leq m$, define $F_i(m)=\sum_{j=0}^{m}F_i(m,j)$ with
$i=1,2,3$, where
\begin{eqnarray*}
F_1(m,j) \hskip-.22cm &=&\hskip-.22cm (-1)^{j}\frac{(2j+2)^2}{(2m+2)^2}\binom{2m+2}{m-j}^2\binom{2m+1}{m}^{-1}, \\
F_2(m,j) \hskip-.22cm &=&\hskip-.22cm (-1)^{j}\frac{(2j+2)^2}{(2m+2)(2m+3)}\binom{2m+3}{m-j}\binom{2m+3}{m-j+1}\binom{2m+2}{m+1}^{-1}, \\
F_3(m,j) \hskip-.22cm &=&\hskip-.22cm
(-1)^{j}\frac{(2j+1)}{(2m+1)}\binom{2m+1}{m-j}^2\binom{2m}{m}^{-1}.
\end{eqnarray*}
Then the creative telescoping algorithm quickly finds the recurrence
\begin{eqnarray*}
F_i(m+1,j)-F_i(m,j) \hskip-.22cm &=&\hskip-.22cm
G_i(m,j+1)-G_i(m,j),
\end{eqnarray*}
where $G_i(m,j)=R_i(m,j)F_i(m,j)$ for $i=1, 2, 3$ and
\begin{eqnarray*}
R_1(m,j) \hskip-.22cm &=&\hskip-.22cm \frac{j^2(j+1)-j(3m+5)(m+1)}{2(j+1)(m-j+1)^2}, \\
R_2(m,j) \hskip-.22cm &=&\hskip-.22cm \frac{j^2(j+1)-j(3m+5)(m+2)}{2(j+1)(m-j+1)(m-j+2)}, \\
R_3(m,j) \hskip-.22cm &=&\hskip-.22cm
\frac{j^3-3j(m+1)^2}{(2j+1)(m-j+1)^2}.
\end{eqnarray*}
If we sum the above recurrence over $j$, we find that the sums
satisfy $F_i(m+1)-F_i(m)=0$ with $i=1, 2, 3$, which easily yield
$F_i(m)=F_i(0)=1$. Hence, (\ref{eqn 4.4}), (\ref{eqn 4.5}) and
(\ref{eqn 4.6}) are proved. \qed

\begin{theorem}
Let $\mathcal{M}=(M_{n,k}(x,y))_{n\geq k\geq 0}$ be given in Section
2. Then for $j=1$ or $2$, there holds
\begin{eqnarray*}
\sum_{k=0}^{n}(-1)^{k}\det\left(\begin{array}{cc}
M_{n,k}(2,2) & M_{n,k+1}(2,2) \\[5pt]
M_{n+j,k}(2,2) & M_{n+j,k+1}(2,2)
\end{array}\right)
\hskip-.22cm &=&\hskip-.22cm 4^{j-1}C_{n+1},
\end{eqnarray*}
or equivalently, for $j=1$ there has
\begin{eqnarray}\label{eqn 4.7}
C_{n+1} \hskip-.22cm &=&\hskip-.22cm
\sum_{k=0}^{n}(-1)^k\frac{(2k+2)(2k+3)(2k+4)}{(2n+2)(2n+3)(2n+4)^2}\binom{2n+4}{n-k}\binom{2n+4}{n-k+1},
\end{eqnarray}
and for $j=2$ there has
\begin{eqnarray}\label{eqn 4.8}
2C_{n+1} \hskip-.22cm &=&\hskip-.22cm \sum_{k=0}^{n}(-1)^k
\frac{(2k+2)(2k+3)(2k+4)}{(2n+2)(2n+3)(2n+6)(2n+7)}\binom{2n+3}{n-k}\binom{2n+7}{n-k+2}.
\end{eqnarray}
\end{theorem}
\pf For $0\leq k\leq n$, define $F_i(n)=\sum_{k=0}^{n}F_i(n,k)$ with
$i=4,5$, where
\begin{eqnarray*}
F_4(n,k) \hskip-.22cm &=&\hskip-.22cm (-1)^k\frac{(2k+2)(2k+3)(2k+4)}{(2n+2)(2n+3)(2n+4)^2C_{n+1}}\binom{2n+4}{n-k}\binom{2n+4}{n-k+1}, \\
F_5(n,k) \hskip-.22cm &=&\hskip-.22cm (-1)^k
\frac{(2k+2)(2k+3)(2k+4)}{2(2n+2)(2n+3)(2n+6)(2n+7)C_{n+1}}\binom{2n+3}{n-k}\binom{2n+7}{n-k+2}.
\end{eqnarray*}
Again the creative telescoping algorithm quickly finds the
recurrence
\begin{eqnarray*}
F_i(n+1,k)-F_i(n,k) \hskip-.22cm &=&\hskip-.22cm
G_i(n,k+1)-G_i(n,k),
\end{eqnarray*}
where $G_i(n,k)=R_i(n,k)F_i(n,k)$ for $i=4, 5$ and
\begin{eqnarray*}
R_4(n,k) \hskip-.22cm &=&\hskip-.22cm \frac{k(k+1)^2-k(3n+7)(n+2)}{(2k+3)(n-k+1)(n-k+2)}, \\
R_5(n,k) \hskip-.22cm &=&\hskip-.22cm
\frac{2k^2(k+2)-2k(3n+7)(n+3)}{(2k+3)(n-k+1)(n-k+3)}.
\end{eqnarray*}
Similarly, if summing the above recurrence over $k$, we find again
that the sums satisfy $F_i(n+1)-F_i(n)=0$ with $i=4, 5$, which
easily generate $F_i(n)=F_i(0)=1$. Hence, (\ref{eqn 4.7}) and
(\ref{eqn 4.8}) are proved. \qed

\begin{remark}
Despite (\ref{eqn 4.3}), (\ref{eqn 4.2C})-(\ref{eqn 4.8}) are proved
by algebraic methods, but bijective proofs are naturally expected,
together with a proof of the following identities involving Motzkin
numbers, namely, for $j=1, 2$, there holds
\begin{eqnarray*}
\sum_{k=0}^{n}(-1)^{n-k}\det\left(\begin{array}{cc}
M_{n,k}(1,1) & M_{n,k+1}(1,1) \\[5pt]
M_{n+j,k}(1,1) & M_{n+j,k+1}(1,1)
\end{array}\right)
\hskip-.22cm &=&\hskip-.22cm 2^{j-1}M_{n}.
\end{eqnarray*}
Besides, another direction in which this research could be taken is
to consider the q-analog of all the identities obtained in this
paper.
\end{remark}

\vskip.2cm

\section*{Acknowledgements} { The work was partially supported by The National Science Foundation of
China and by the Fundamental Research Funds for the Central
Universities.}

\vskip.2cm


\end{document}